\documentclass[12pt]{article}

\usepackage{graphicx}
\usepackage{psfrag}
\usepackage{epsf}
\usepackage{amsmath,amsfonts,amssymb,latexsym,amsthm}
\usepackage{enumitem}
\usepackage[notref,notcite]{showkeys}
\usepackage{algorithmic}
\usepackage{algorithm}
\usepackage[width=160mm,height=240mm,left=35mm,foot=10mm]{geometry}
\usepackage{setspace}
\newcommand{\ignore}[1]{}
\newcommand{\startClaims}{\setcounter{claim}{0}}
\newtheorem{theorem}{Theorem}[section]

\newtheorem{lemma}[theorem]{Lemma}

\newtheorem{conjecture}[theorem]{Conjecture}

\usepackage[usenames]{color}
\usepackage{multirow}

\newcommand{\calD}{\mathcal{D}}
\newcommand{\calP}{\mathcal{P}}
\newcommand{\calS}{\mathcal{S}}

\title{On the crossing number of $K_{13}$}

\author{ \textsc{Dan McQuillan, Shengjun Pan, and R.~Bruce Richter}\ignore{\footnotemark[2]}\\[0.25em]
{\small{Department of Mathematics}} \\[-0.25em]
     {\small{Norwich University}} \\[-0.25em]
     {\small{Vermont}} \\[-0.1em]
     {\small\texttt{dmcquill@norwich.edu}}\\
     {\small{San Diego}}\\[-0.25em]
     {\small{California}}\\[-.01em]
     {\small\texttt{shengjunpan@gmail.com}}\\
     {\small{Department of Combinatorics \& Optimization}} \\[-0.25em]
     {\small{University of Waterloo}} \\[-0.25em]
     {\small{Canada}} \\[-0.1em]
     {\small\texttt{brichter@uwaterloo.ca}}}

\date{\today}

\newcommand{\eopf}{\raisebox{0.8ex}{\framebox{}}}

{\noindent{\bf Sketch of a proof.}\ }%
{\hfill\eopf\par\bigskip}%
{\noindent{\bf Ideas for a proof.}\ }%
{\hfill\eopf\par\bigskip}%
{\noindent{\bf Proof in progress.}\ }%
{\hfill\eopf\par\bigskip}%

\newenvironment{cproof}
{\noindent{\bf Proof.}\startClaims\ }
{\hfill{{\rule[0.0ex]
{7pt}{7pt}}}\par\bigskip}

\def\crn{{\hbox{\rm cr}}}

\begin{document}

\maketitle

\begin{abstract}  Since the crossing number of $K_{12}$ is now known to be 150, it is well-known that simple counting arguments and Kleitman's parity theorem for the crossing number of $K_{2n+1}$ combine with a specific drawing of $K_{13}$ to show that the crossing number of $K_{13}$ is one of the numbers in $\{217,219,221,223,225\}$.  We show that the crossing number is not 217.  \end{abstract}

\section{Introduction and drawings of $K_{12}$ having $151$ crossings}\label{intro}

The {\em crossing number\/} $\crn(D)$ of a drawing $D$ of a graph $G$ is the number of pairwise intersections of edges of $G$ in $D$.  The {\em crossing number\/} $\crn(G)$ of $G$ is the least $\crn(D)$ over all drawings $D$ of $G$ in the plane.   A long-standing conjecture is that the crossing number of the complete graph $K_n$ satisfies the following formula.

\begin{conjecture}\label{co:zaran}  For $n\ge 1$, 
$$\crn(K_n)=\frac14 \left\lfloor {\mathstrut\frac { \mathstrut{n}}2}\right\rfloor \left\lfloor \frac {n-1}2\right\rfloor \left\lfloor \frac {n-2}2\right\rfloor \left\lfloor \frac {n-3}2\right\rfloor\,.$$
\end{conjecture}

Recently, Pan and Richter \cite{pr} used a computer to prove that the crossing number of $K_{11}$ is indeed 100.  A simple counting argument shows that $\crn(K_{12})=150$ (as long as the conjecture holds for $K_{2n-1}$, it automatically holds for $K_{2n}$).   The next case is $K_{13}$.

We define $Z(n)$ by the formula
$$
\frac14 \left\lfloor {\mathstrut\frac { \mathstrut{n}}2}\right\rfloor \left\lfloor \frac {n-1}2\right\rfloor \left\lfloor \frac {n-2}2\right\rfloor \left\lfloor \frac {n-3}2\right\rfloor\,.
$$
There is a standard drawing that shows $\crn(K_n)\le Z(n)$.  

Let $D$ be a drawing of $K_n$.  For each vertex $v$ of $K_n$, $D-v$ is the drawing of $K_{n-1}$ obtained from $D$ by removing $v$ and its incident edges.  Since each crossing of $D$ is in $n-4$ of the $D-v$, we see that

\begin{eqnarray}\label{eq:counting} (n-4)\,\crn(D) = \sum_{v\in V(K_n)}\crn(D-v)\,.\end{eqnarray}
In particular, $(n-4)\crn(K_n)\ge n\hskip 1pt\crn(K_{n-1})$.

Kleitman \cite{kleitman} proved that any two drawings of $K_{2n+1}$ with no tangencies and no two edges both incident with a common vertex and crossing each other have crossing numbers having the same parity.  Together with Equation (\ref{eq:counting}), the crossing number of $K_{13}$ is one of 217, 219, 221, 223, and 225.  

It is the purpose of this work to show that $\crn(K_{13})>217$, the first substantial progress on this number since it was shown that $\crn(K_{13})\le 225$.

A second simple idea is that of {\em duplicating a vertex\/} in a drawing of $K_{n}$ to obtain a drawing of $K_{n+1}$.  Let $v$ be any vertex in a drawing $D$ of $K_{n}$.  We denote by $\crn(v,D)$ the number of crossings in $D$ involving at least one edge incident with $v$.   By placing a new vertex $v'$ near $v$ in the drawing and then making the edges from $v'$ parallel to those of $v$, we obtain a drawing $D'$ of $K_{n+1}$.  The following inequality, which we state only when $n$ is odd, is now standard (see for example Woodall \cite{woodall} and Christian, Richter, and Salazar \cite{crs}):

\begin{eqnarray}\label{eq:dupVert}  \crn(D')\ge \crn(D)+\crn(v,D)+2\binom{\frac {n-1}2}2\,.\end{eqnarray}
Moreover, equality can always be achieved.

A surprisingly interesting parameter for a drawing $D$ of $K_n$ is the quantity $\delta(D) = \crn(D)-Z(n)$.
 A drawing $D$ of $K_{2m}$ has the {\em normal deficiency property\/}  if, for every vertex $v$ of $K_{2m}$,  $\delta(D-v)\le 2\delta(D)$. We will usually   abbreviate this property to NDP.  For example, every optimal drawing of $K_{12}$ has NDP, using Equation (\ref{eq:counting}) to show that each $K_{11}$ in $D$ has 100 crossings.

One point about the NDP property is that if there is a drawing $D$ of $K_{2n}$ and a vertex $w$ so that $\delta(D-w)>2\delta(D)$, then $\crn(K_{2n+1})<Z(2n+1)$.  This is because duplicating $w$ in $D$ yields a drawing $D'$ of $K_{2n+1}$ with $\crn(D')<Z(2n+1)$.      In particular, if there is a drawing $D$ of $K_{12}$ that does not have NDP, then $\crn(K_{13})<225$.   The following is a converse if $\crn(K_{13})=217$.  

\begin{lemma}\label{lm:theory} Suppose $D$ is a drawing of $K_{13}$ with 217 crossings.  Then:
\begin{enumerate}
\item\label{it:three} there are 3 vertices $v$ of $K_{13}$ so that $\delta(D-v)=1$;
\item\label{it:ten} there are 10 vertices $v$ of $K_{13}$ so that $\delta(D-v)=0$; and
\item\label{it:notNDP} if $u$ and $v$ are distinct vertices so that $\delta(D-u)=\delta(D-v)=1$, then $\delta(D-u-v)=4$.
\end{enumerate}
\end{lemma}

\begin{cproof}  We know that $\crn(K_{12})=Z(12)=150$.  Equation \ref{eq:counting} shows $\sum_{v\in V(K_{13})}\crn(D-v)=9\hskip 1pt\crn(D)$.  Therefore, $9(217) = a(150)+b(151)+c(152)+\cdots$.  From this equation, it follows that $a\ge 10$ and $b+2c+3d=3$.    

It follows that there is a $v\in V(K_{13})$ so that $\crn(D-v)>150$ and so, from (\ref{eq:counting}), there is a $u\in V(K_{13})$, with $u\ne v$, so that $\crn(D-v-u)>100$.  Going the other way, we deduce that $\crn(D-u)>150$.

Of the ${13\choose 11}$ $K_{11}$'s in $D$, we have found at least one that has more than 100 crossings.  
If there is only the one $K_{11}$ with more than 100 crossings, then  $\crn(D-w)=150$ for all 11 vertices $w$ of $K_{13}$ other than $u$ and $v$.  It follows that all the other ${13 \choose 11}-1$ $K_{11}$'s in $D$ have exactly 100 crossings.  But this implies $8\hskip1pt\crn(D-u)=\crn(D-u-v)+11\cdot 100= 8\hskip1pt\crn(D-v)$.  In particular, $\crn(D-u)=\crn(D-v)$.  Therefore, one of $b$, $c$, and $d$ is 2, and yet $b+2c+3d=3$, a contradiction.  

It follows that there is another $K_{11}$ with more than 100 crossings.  Therefore, there is a third vertex $w$ of $K_{13}$ so that $\crn(D-w)>150$ as well.  In this case, $b+c+d\ge 3$.  Since $b+2c+3d=3$, we have $b=3$, proving (\ref{it:three}) and (\ref{it:ten}).

{\bf\large (FIX in REVISION)} that completes the proofs of (\ref{it:three}) and (\ref{it:ten}).

For (\ref{it:notNDP}), let $v_1$, $v_2$, and $v_3$ be the three vertices $v$ of $K_{13}$ so that $\crn(D-v)=151$.   If $\{i,j,k\}=\{1,2,3\}$, then 
$$
8\crn(D-v_i)=\crn(D-v_i-v_j)+\crn(D-v_i-v_k) + 10\cdot 100\,.
$$
We conclude that $\crn(D-v_i-v_j)+\crn(D-v_i-v_k)=208$.  Solving these three equations yields $\crn(D-v_i-v_j) =104$, for each distinct $i,j\in \{1,2,3\}$.  
\end{cproof}

Lemma \ref{lm:theory} shows that, if $\crn(K_{13})=217$, then there is a drawing $D$ of $K_{12}$ having 151 crossings so that some $K_{11}$ has 104 crossings in $D$.  Conversely, if such a drawing of $K_{12}$ exists, then duplicating a vertex shows $\crn(K_{13})<225$.  The remainder of this work is explaining how we show that no such $D$ exists.

\section{Equivalent drawings of $K_n$}

In this section, we describe our algorithm for determining that no drawing of $K_{12}$ having 151 crossings contains a $K_{11}$ having 104 crossings.

It is well-known that some drawing $D$ of a graph $G$ with $\crn(D)=\crn(G)$ is {\em good\/}, which means the following properties hold:
\begin{enumerate}
\item no two edges incident with a common vertex cross in $D$;
\item no three edges are concurrent at a point; and 
\item no two edges intersect more than once.
\end{enumerate}

We will check every good drawing of $K_{12}$ having 151 crossings to see if it contains a subdrawing of $K_{11}$ having 104 crossings.  If there is one, then $\crn(K_{13})<225$, while if there is not, then $\crn(K_{13})>217$. 

Equation \ref{eq:counting} shows that, if $D$ is a drawing of $K_{12}$ having 151 crossings, then at least one of the twelve subdrawings of $K_{11}$ must have precisely 100 crossings.  Repeating this argument with a drawing $D'$ of $K_{11}$ having 100 crossings, one of the eleven subdrawings of $K_{10}$ must have at most 63 crossings.

Thus, our program is simple in conception.  Generate all drawings of $K_{10}$ having at most 63 crossings.  For each, try to extend it to an optimal $K_{11}$ and then to a $K_{12}$ with 151 crossings.   Then check each of the twelve vertex-deleted subdrawings of $K_{12}$ to see if any has 104 crossings.

This program is analogous to the one employed by Pan and Richter \cite{pr} to show that $\crn(K_{11})\ge 100$.  However, in that instance, it was not necessary to go beyond testing, for each drawing of $K_{10}$ having at most 62 crossings, each face to see if an eleventh vertex could be placed there to get a drawing of $K_{11}$ having fewer than 100 crossings.  In this context, we need to keep the $K_{11}$'s to create the $K_{12}$'s; in particular, the number of optimal drawings of $K_{11}$ is still not known.  Such a direct approach did not succeed for us, as too many drawings of $K_{11}$ were produced.  We were able to overcome this obstacle by considering a single representative of a set of equivalent drawings.

Two good drawings $D$ and $D'$ of $K_n$ are {\em equivalent\/} if the pairs of crossing edges are the same in both $D$ and $D'$.  For us, this comes up in the context of extending a drawing $D_0$ of $K_{10}$ in all possible ways to optimal drawings of $K_{11}$ and then to drawings of $K_{12}$ having 151 crossings.   

In the step from $K_{10}$ to $K_{11}$, two drawings $D$ and $D'$ differ only in how the edges incident with the new vertex $v$ are drawn.  Consider the drawing of an edge $D'(vw)$ relative to $D(K_{11})$.  Equivalence means that there is  one region of $D(vw)\cup D'(vw)$ that has all the other vertices of $D(K_{10})$, while the remaining regions bound only crossings of $D(K_{10})$, and that this is true for all ten edges incident with $v$.

In the step from $K_{11}$ to $K_{12}$, we have the same discussion.  However, there is another consideration to take into account.  Two drawings of $K_{11}$ that are equivalent relative to the $K_{10}$ may become inequivalent once the twelfth vertex is added, as this new vertex may lie in a region of $D(vw)\cup D'(vw)$ different from the one containing the other 9 vertices of $K_{10}$.

In our program, we start with the drawing $D_0$ of $K_{10}$.  Dijkstra's algorithm (that is, Breadth First Search) is used to determine, for each face $F$ of $D_0$, whether the eleventh vertex $v$ can be placed there to obtain an optimal drawing of $K_{11}$.   

For each edge $vw$, we find all optimal routings and these are sorted into equivalence classes, determined by which edges of $D_0$ are crossed in the routing.  The optimal drawings of $K_{11}$ are obtained by choosing a routing for each of the ten edges $vw$.  For each edge $vw$ and each equivalence class $[d(vw)]$ of drawings of $vw$, we choose a representative drawing of $vw$; these 10 representatives are then combined to yield a drawing of $K_{11}$.  Any optimal drawing of $K_{11}$ is equivalent to one of these.

There is one subtlety here that needs discussion.  The selection of the representatives may result in routings of the ten edges incident with $v$ that are not mutually compatible; the routings may force two edges incident with $v$ to cross.  This is something our program checks for, but fortunately never came up and so we never had to deal with it.  Thus, we have not been required to determine whether it is possible that one pair of representatives for particular equivalence classes for two edges might be tangled, while a different pair of representatives is not tangled.  The representatives we happened to choose were always untangled, so this problem did not arise in the execution of the program.

We now deal with the possibility of equivalent (relative to $D_0$ and $F$) drawings of $K_{11}$ becoming inequivalent with the introduction of the twelfth vertex $v'$.  Let $F'$ be the face of $D_0$ into which $v'$ is placed.  If two $D_0$-equivalent routings of an edge $vw$ incident with $v$ are now inequivalent, it is because they separate $v'$ from the other nine vertices of $K_{10}$.  This implies that they have a $K_{10}$-equivalent routing {\em that goes through $F'$\/}.  

For each edge $vw$ for which it is possible, we choose a $K_{10}$-equivalent routing through $F'$ to obtain our optimal drawing of $K_{11}$.   These routings will break up $F'$ into at most 10 sub-faces and we simply try placing $v'$ into each of these to make our drawing of $K_{12}$.  In particular, if $F'=F$, then we will have 10 sub-faces to try.    These different possibilities for $v'$ include representatives of all extensions of equivalents of this $K_{11}$ to a $K_{12}$ in which $v'$ is placed in $F'$.  We formalize this important point as a theorem.

\begin{theorem}  Let $D$ be a drawing of $K_{12}$ having 151 crossings.  Let $v_{11}$ and $v_{12}$ be vertices of $K_{12}$ so that, among all $K_{10}$'s contained in $D$,  $D-\{v_{11},v_{12}\}$ has fewest crossings so that $\crn(D-v_{12})=100$.  Then some drawing equivalent to $D$ is produced by our algorithm from an isomorph of $D-\{v_{11},v_{12}\}$.  \end{theorem}

\begin{cproof}  We have already seen that the drawing $D-\{v_{11},v_{12}\}$ has at most 63 crossings and so is one of the drawings of $K_{10}$ that we consider.  Let $F_{11}$ and $F_{12}$ be the faces of $D-\{v_{11},v_{12}\}$ containing $v_{11}$ and $v_{12}$, respectively.

For each edge $v_{11}w$, with $w$ in $K_{12}-\{v_{11},v_{12}\}$, we check to see if there is a $D[K_{12}]$-equivalent of $v_{11}w$ routing through $F_{12}$.  One at at time, we will shift such edges, with the end result being a good drawing of $K_{12}$, equivalent to $D$, with all possible edges incident with $v_{11}$ drawn through $F_{12}$.    The process is an induction on the number of edges incident with $v_{11}$ that can be equivalently routed (relative to $K_{12}$) through $F_{12}$ but are not so in the current drawing.  The initial current drawing is $D$.

Let $v_{11}w$ be an edge that can be equivalently routed through $F_{12}$, but is not so routed in our current drawing $D^*$.  Let $e$ be an equivalent routing through $F_{12}$.  If simply replacing $D^*[v_{11}w]$ with $e$ results in a good drawing of $K_{12}$, then we do this.  

Otherwise, there is an edge $ab$ that crosses $e$ in violation of goodness.  Then $ab$ and $e$ bound a digon whose interior contains only crossings.  Choosing $ab$ to bound a minimal such digon (in the sense of no other digon is contained inside it), we can equivalently shift the portion of $ab$ in the boundary of this minimal digon to the other side of the portion of $e$ in the boundary of this digon.  This reduces the total violation of goodness, so repetition results in a new, equivalent, good drawing $D'$ of $K_{12}$ that uses $e$ in place of $v_{11}w$.

This new drawing $D'$ has fewer edges incident with $v_{11}$ that can be, but are not, routed through $F_{12}$ than $D^*$ has.  By induction, there is an equivalent drawing with all possible edges routed through $F_{12}$.  This drawing is equivalent to one we consider.
\end{cproof}

In going from $K_{11}$ to $K_{12}$, there is one additional complication.  We are looking for drawings of $K_{12}$ having 151 crossings.  This is one more than optimal.  For each of the faces $F''\subseteq F'$ of drawings of $K_{11}$, we again use Dijkstra's algorithm to determine the minimum number of additional crossings added to obtain the drawing of $K_{12}$.  Obviously, if this minimum number yields a drawing of $K_{12}$ with more than 151 crossings, we ignore this face $F''$.  If the minimum number produces a drawing with precisely 151 crossings, we again check all equivalence classes for the edges and see whether some subdrawing of $K_{11}$ has precisely 104 crossings.

In the final case, the face $F''$ produces a drawing of $K_{12}$ having 150 crossings.  In such a case, we choose an edge $v'w$ incident with $v'$ (eventually testing all 11).  We then test whether there is a rerouting of $v'w$ to have one additional crossing, and find all such routings, divided into equivalence classes.  Finally, we check as above, whether, given a selection from each equivalence class, there is a subdrawing of $K_{11}$ having 104 crossings.

We executed this program and the outcome was that, among all the considered drawings of $K_{12}$ having 151 crossings, none contains a $K_{11}$ having 104 crossings.  As indicated, this implies that no drawing of $K_{12}$ having 151 crossings has a subdrawing of $K_{11}$ having 104 crossings.  Lemma \ref{lm:theory} implies that there is no drawing of $K_{13}$ having 217 crossings, and, therefore, $\crn(K_{13})>217$.

%========= ADDED BY SHENGJUN ========>

\section{Generating drawings of $K_{n+1}$ from drawings of $K_n$}
For any integer $n\ge 3$, let $\calD_n^{c}$ be the set of all good
drawings of $K_n$ that have exactly $c$ crossings. Similarly, let
$\calD_n^{\le c}$ ($\calD_n^{\ge c}$) be the set of all good drawings
of $K_n$ that have at most (at least) $c$ crossings. The following results are immediate from Equation \ref{eq:counting}.
\begin{theorem}
\label{thm:contains}
Let $\calS_4,\calS_5,\dotsc,\calS_{12}$ be the sets of drawings:
\[
\calD_4^0,\,\,
\calD_5^1,\,\,
\calD_6^3,\,\,
\calD_7^9,\,\,
\calD_8^{\le 20},\,\,
\calD_9^{\le 36},\,\,
\calD_{10}^{\le 63},\,\,
\calD_{11}^{\le 100},\,\,
\calD_{12}^{151},
\]
respectively. Then, for each $n=4,5,\dotsc,11$, for any drawing $D$ in
$\calS_{n+1}$, there is a vertex $v$ of $K_{n+1}$ so that $D-v\in
\calS_n$. \end{theorem}

It was shown in~\cite{pr} that $\crn(K_{11})=100$. Thus $\calD_{11}^{\le
  100} = \calD_{11}^{100}$, although it is not required in this
paper. In fact, this paper re-proves the result $\crn(K_{11})=100$ by showing
that $\calD_{11}^{\le 100} = \calD_{11}^{100}$ .

Recall that our objective is to generate all drawings in
$\calD_{12}^{151}$ and show that none of them has a subdrawing in
$\calD_{11}^{\ge 104}$. Theorem~\ref{thm:contains} shows how to 
generate all drawings in $\calD_{12}^{151}$ from $\calD_4^0$ by
repeatedly adding new vertices and new edges.

To state the procedure more precisely, given positive integers $n$ and $c$ we introduce a generic
algorithm that takes a set $\calD_n$ of drawings of $K_n$, and outputs
all drawings of $K_{n+1}$ that have at most $c$ crossings and a
subdrawing in $\calD_n$. The steps are
listed in Algorithm~\ref{alg:insertion}, where the distance $d(F,v_i)$
is the smallest length of any routing from $F$ to $v_i$.

We supply $\calD_4^0$ (consisting of a single drawing) as the input to
Algorithm~\ref{alg:insertion} to obtain $\calD_5^1$, which in turn
becomes the new input to get $\calD_6^3$. Continuing the process, theoretically we
could obtain $\calD_{12}^{151}$.

\begin{algorithm}
\caption{Generating drawings of $K_{n+1}$ from drawings of $K_n$}
\label{alg:insertion}
\algsetup{
    linenodelimiter=.
}
\begin{algorithmic}[1]
\REQUIRE $\calD_n$: a set of drawings of $K_n$
\ENSURE $\calD_{n+1}^{\le c}$: the set of all drawings of $K_{n+1}$ that
\begin{itemize}\setlength{\itemsep}{-1ex}
\item have at most $c$ crossings, and
\item have a subdrawing in $\calD_n$.
\end{itemize}
\STATE Let $\calD_{n+1}^{\le c}\leftarrow\emptyset$
\FORALL{$D\in \calD_n$}
\STATE Let $\varepsilon\leftarrow c - \crn(D_n)$\\
Let $v_1,v_2,\dotsc,v_n$ be the vertices of $D$
\FORALL{face $F$ in $D$} \label{step:face}
\FORALL{$i=1,2,\dotsc,n$}
\STATE Compute distances $d(F,v_i)$.
\STATE Find set $\calP_i$ of all routings of length at most
$d(F,v_i)+\varepsilon$ from $F$ to $v_i$. \label{step:routing}
\ENDFOR
\FORALL{products of routings $(P_1,P_2,\dotsc,P_n)\in \calP_1\times \calP_2\times\dotsb \times \calP_n$}
\IF {the total length is at most $c - \crn(D)$, and the $P_j$ do not cross each other}
\STATE Insert a new vertex $v_{n+1}$ in $F$.
\STATE Generate a new drawing $D'$ by drawing $P_1,P_2,\dotsc,P_n$ from $v_{n+1}$ to
$v_1,v_2,\dotsc,v_n$ respectively.
\IF {$D'$ is not isomorphic to any drawing already in the output} \label{step:iso}
\STATE  Add $D'$ to $\calD_{n+1}^{\le c}$.
\ENDIF
\ENDIF
\ENDFOR
\ENDFOR
\ENDFOR
\end{algorithmic}
\end{algorithm}
\subsection{Generating drawings of $K_{11}$ from drawings of $K_{10}$}
\label{subsec:10to11}
In fact, however, we only use 
Algorithm~\ref{alg:insertion} to obtain $\calD_{10}^{\le
  63}$.  In order to get from $\calD_n^{\le c}$ to $\calD_{n+1}^{\le c'}$, we separately do each $\calD_n^{c''}$, for each $c''\le c$ (and $\ge \crn(K_n)$).  Table~\ref{tbl:count_and_cost} lists the sizes of the sets and
the computational time to produce them.

\begin{table}[h]
  \[\renewcommand{\arraystretch}{1.2} \begin{array}{|c|r|l|} \hline
    \text{drawings} & \text{\# drawings} & \text{cost of time}\\
    \hline \calD_4^0 & 1 & \\ \hline
    \calD_5^1 & 1 & \multirow{3}{*}{$\le 1$ second}\\ \cline{1-2}
    \calD_6^3 & 1 & \\  \cline{1-2}
    \calD_7^9 & 1 & \\  \hline
    \calD_8^{18} & 3 & \multirow{3}{*}{5 seconds}\\
    \calD_8^{19} & 18 & \\
    \calD_8^{20} & 88 & \\ \hline \calD_{9}^{36} & 3,080 & 7 \text{
      minutes}\\ \hline
    \calD_{10}^{60} & 5,679 & \multirow{4}{*}{505 hours}\\
    \calD_{10}^{61} & 115,095 & \\
    \calD_{10}^{62} & 1,080,968 & \\
    \calD_{10}^{63} & 6,171,344 & \\ \hline \end{array} \]
  \caption{Number of drawings and time costs}
  \label{tbl:count_and_cost} \end{table} 

The number of drawings in $\calD_{11}^{100}$ is too large for the resources available to us; we were not able to generate $\calD_{11}^{100}$. For example, Table~\ref{tbl:test_K11} shows the
number of drawings in $\calD_{11}^{100}$, and corresponding 
computation time, generated by Algorithm~\ref{alg:insertion} with
input of \emph{one} drawing from $\calD_{10}^{\le 63}$. The test is
repeated four times with the input taken from
$\calD_{10}^{60},\calD_{10}^{61},\calD_{10}^{62}$ and
$\calD_{10}^{63}$, respectively. 

\begin{table} 
\[
  \begin{array}{|c|c|c|c|c|c|}
  \hline
\multicolumn{2}{|c|}{\text{input}}
 & \calD_{10}^{60} &\calD_{10}^{61} &\calD_{10}^{62} & \calD_{10}^{63} \\ \hline
\multirow{2}{*}{Algorithm~\ref{alg:insertion}}
& \text{number of $K_{11}$'s} & $653,125$&$310,150$&$73,261$ & $2,147$\\ \cline{2-6}
& \text{running time} &60.3 \text{ hours}&16.3 \text{ hours} &7.4 \text{ hours} & 2.0 \text{ hours} \\ \hline
\hline
\multirow{2}{*}{Algorithm~\ref{alg:compound}}
& \text{number of $K_{11}$'s} &$69,064$&$36,946$&$12,771$&$1,400$\\ \cline{2-6}
%& \text{\# non-ismorphic} &6,165&3,736&1,204&113\\ \cline{2-6}
& \text{running time} &1.5\text{ minutes}&1.5 \text{ minutes}& 40 \text{ seconds}&30 \text{ seconds}\\ \hline
\end{array}
\]
\caption{Tests on numbers of $K_{11}$'s from a single drawing of $K_{10}$}
\medskip
%\centerline{\small \emph{Note:} Each test takes a single drawing of $K_{10}$ as input.}
\label{tbl:test_K11}
\end{table}

%\subsubsection{Equivalent routings}
As discussed in Section 2, in order to prove that $\crn(K_{13})>217$, it suffices to consider only inequivalent drawings of $K_{11}$ having 100 crossings.  Thus, we shall generate a set 
  $\widetilde\calD_{11}^{100}$ of drawings of $K_{11}$ having 100 crossings and containing at least one representative from each equivalence class of drawings in $\calD_{11}^{100}$.   Algorithm~\ref{alg:compound}
  gives the exact steps. It is a
  modified version of Steps~\ref{step:partition} --
  \ref{step:choose_end} of Algorithm~\ref{alg:insertion}.  These modifications are essential  in reducing the number
  of drawings of $K_{11}$ that we have to check.

  As shown in Step~\ref{step:entangled}, if the
  routings chosen from the previous step are entangled, we simply save
  the current drawing of $K_{10}$ to $\calD_{10}^{\rm error}$. 
Since we need to produce good drawings, our analysis does not apply for the drawings in $\calD_{10}^{\rm error}$.  For a drawing in $\calD_{10}^{\rm error}$, we would
  consider all products of routings, equivalent or non-equivalent, as
  long as they are not entangled. If the size of $\calD_{10}^{\rm
    error}$ were large, we might not be able to complete the computation. Fortunately, our
  program never selected entangled routings and, therefore, at the end of the computation,  $\calD_{10}^{\rm error}=\emptyset$, so no additional work is required.

\begin{algorithm}
\caption{Generating drawings in $\widetilde\calD_{11}^{100}$ from drawings of $\calD_{10}^{\le63}$}
\label{alg:compound}
\algsetup{
    linenodelimiter=.
}
\begin{algorithmic}[1]
\REQUIRE $\calD_{10}^{\le 63}$: the set of drawings of $K_{10}$ with at
most 63 crossings
\ENSURE \ 
\begin{itemize}\setlength{\itemsep}{0pt}
\item $\widetilde\calD_{11}^{100}$: drawings of $K_{11}$ with 100
crossings obtained by adding only non-equivalent routings.
\item $\calD_{10}^{\rm error}$: a subset of drawings of $\calD_{10}^{\le
    63}$ that would have entangled routings.
\end{itemize}
\STATE Let $\widetilde\calD_{11}^{100}\leftarrow \emptyset$\\
Let $\calD_{10}^{\rm error}\leftarrow\emptyset$
\FORALL{ $D\in \calD_{10}^{\le 63}$}
\STATE Let $\varepsilon\leftarrow 151 - \crn(D)$
\STATE Let $v_1,v_2,\dotsc,v_{10}$ be the vertices of $D$
\FORALL{faces $F$ of $D$}\label{step:face2}
\FORALL{$i=1,2,\dotsc,10$}
\STATE Compute distances $d(F,v_i)$.
\STATE Find set $\calP_i$ of all routings of length at most
$d(F,v_i)+\varepsilon$ from $F$ to $v_i$. \label{step:routing2}
\STATE \label{step:partition} Partition $\calP_i$ into equivalence classes:
$\calP_i=\bigcup_{j} \calP_{i,j}$
\ENDFOR
\FORALL{products of classes  $(\calP_{1,j_1},\calP_{2,j_2},\dotsc,\calP_{10,j_{10}})$}
\STATE \label{step:choose_begin} Choose any $P_i\in \calP_{i,j_i}$ as follows:
\IF{there exists a routing in $\calP_{i,j_i}$ passing through $F$}
\STATE let $P_i$ be any such routing.
\ELSE
\STATE Choose an arbitrary routing in $\calP_{i,j_i}$ for $P_i$.
\ENDIF  \label{step:choose_end} 
\IF {the total length of all routings is at most $ c - \crn(D)$}
\IF {the paths $P_1,P_2,\dotsc,P_{10}$ do not cross one another}\label{step:entangled} 
\STATE insert a new vertex $v_{11}$ in $F$ and obtain $D'$
 by drawing $P_1,P_2,\dotsc,P_{10}$
from $v_{11}$ to $v_1,v_2,\dotsc,v_{10}$, respectively.
%\IF {$D'$ is not isomorphic to any drawing already in the output} \label{step:iso2}
\STATE  Add $D'$ to $\widetilde\calD_{11}^{100}$ \label{step:output}
%\ENDIF
\ELSE
\STATE Add $D$ to $\calD_{10}^{\rm error}$
\ENDIF
\ENDIF
\ENDFOR
\ENDFOR
\ENDFOR
\end{algorithmic}
\end{algorithm}

%\subsubsection{Skipping isomorphism testing} 

In Algorithm~\ref{alg:insertion} a drawing is
tested for isomorphism against previously saved drawings; only
non-ismorphic drawings are saved. This is necessary for producing
drawings of $K_{n}$ for $n$ up to 10; otherwise we would have even
more drawings of $K_{11}$ to consider.

Algorithm~\ref{alg:compound} does not test each drawing of $K_{11}$ against previously saved drawings. Therefore, it may test the same $K_{11}$ isomorph more than once. In spite of this potential duplication, the overall computation time is greatly reduced.  When a new $K_{11}$ is generated from a drawing $D$ of $K_{10}$ having at most 63 crossings, the other 10 vertices are checked in turn to see if their deletion produces a drawing $D'$ of $K_{10}$ with $\crn(D')<\crn(D)$.  If such a $D'$ is found, then the drawing of $K_{11}$ is discarded, as an equivalent drawing must have been produced when considering $D'$ earlier.

Table~\ref{tbl:test_K11} lists the number of drawings in
$\calD_{11}^{100}$ and the corresponding computational time from
Algorithm~\ref{alg:compound}. As we can see, even allowing isomorphic
drawings of $K_{11}$ in the output,  both the size of the output and
the running time are significantly reduced.

\subsection{Generating drawings of $K_{12}$ form drawings of $K_{11}$}

For each drawing  obtained from
Algorithm~\ref{alg:compound}, we would like to generate all possible
drawings in $\calD_{12}^{151}$ by adding a new vertex $v_{12}$ and 11
routings $P_1,P_2,\dotsc, P_{11}$ from $v_{12}$ to existing vertices,
and then check if the resulting new drawing has a subdrawing in
$\calD_{11}^{\ge 104}$.

Equivalently, for each $D\in\widetilde\calD_{11}^{100}$, we would like
to check if there is a drawing in $\calD_{12}^{151}$ having a vertex whose incident edges combine for a total of at most 47 crossings.  These crossings are easily determined from $D$
itself and the crossings of the new routings. Thus, we
do not actually generate any drawings of $K_{12}$.

Furthermore, entanglements of the eleven new edges are irrelevant; 
untangling them does not change the number of crossings at one of the original $K_{11}$, which is what we count.  Thus we may also skip checking for
entanglement in going from $K_{11}$ to $K_{12}$.

In our implementation, the generation of drawings of $K_{12}$ is
integrated into Algorithm~\ref{alg:compound}: at
Step~\ref{step:output} we generate and test drawings of $K_{12}$
from the current drawing $D$ of $K_{11}$ for generating drawings of
$K_{12}$ right away, and hence there is no need to save $D$ for later use.

\section{Implementation}

%\subsection{Setup}
Note that both Algorithms~\ref{alg:insertion} and~\ref{alg:compound}
can be made parallel by distributing input drawings to different
computing resources. Our program is set up to run on the cluster
\texttt{saw.sharcnet.ca} with $2712$ cores of 2.83 GHZ. The actual
number of cores allocated for us was 256, which is sufficient for the
program to finish in a reasonable amount of time (See
Table~\ref{tbl:count_and_cost}). {\bf\large(Revision makes ``See" into ``see".)}

%\subsection{Isomorphism testing}
Our program also utilizes the C programming
package~\emph{Nauty}~\cite{nauty} for graph-isomorphism testing. The
isomorphism testing is used by Step~\ref{step:iso} in
Algorithm~\ref{alg:insertion}. As explained in
Subsection~\ref{subsec:10to11}, we do not use isomorphism testing on
the output of Algorithm~\ref{alg:compound}.

On the other hand, for both Algorithms~\ref{alg:insertion} (Step~\ref{step:face})
and~\ref{alg:compound} (Step~\ref{step:face2}), we use isomorphism
testing to eliminate isomorphic faces.

%\subsection{Optimization}
 In searching for routings
at Step~\ref{step:routing} of Algorithm~\ref{alg:insertion} and
Step~\ref{step:routing2} of Algorithm~\ref{alg:compound}, it is an important saving that we need only consider routings that go through distinct faces. 

In general, a routing from a face $F$ to a vertex $v$ could potentially pass
through the same face $F'$ more than once when its length is greater
than $d(F,v)$. In other words, the routing is a walk in the dual graph
but not a path.   The following is a slight improvement to \cite[Claim 1, p.~132]{pr} that suffices for our purposes.

\begin{theorem}  Let $n\ge 5$ be an integer, let $D$ be a good drawing of $K_{n+1}$, let $v$ be any vertex of $K_{n+1}$ and let $D_n=D-v$.   For a vertex $w$ of $K_{n+1}-v$,  let $d(w)$ denote the fewest possible crossings in a routing in $D_n$ from the face containing $v$ to $w$. 
\begin{enumerate}\item
Then  any edge $vw$ that has fewer than $d(w) + n-1$ crossings in $D$ goes through only distinct faces of $D_n$. \item    In particular, if the number of crossings of $vw$ in $D$ is at most $d(w)+3$, then the routing of $vw$ in $D$ goes through distinct faces of $D_n$.
\end{enumerate}
\end{theorem}

\begin{cproof}  Suppose $(F_0,F_1,\dots,F_{t-1},F_0)$ is a cyclic sequence of distinct faces $F_0,F_1,\dots,F_{t-1}$ that occurs consecutively in the routing  $D[vw]$.  Let $C$ be a simple closed curve representing this sequence, in the sense that $C$ crosses precisely the edge segments that are coincident with $F_{i-1}$ and $F_i$, for each $i=1,2,\dots,t$, where the indices are read modulo $t$.

Let $n_1$ and $n_2$ be the numbers of
vertices on the two sides of $C$; we note $n_1\ge0$, $n_2\ge0$ and
$n_1+n_2=n\ge 5$.   If, say $n_1=0$, then every edge of $K_n$ crosses $C$ an 
even number of times (into and out from the side with no vertices).  As $D$ is a good drawing,  $D[vw]$ crosses each edge of $D_n$ at most once.  Therefore, $C$ crosses no edge of $D_n$, showing the sequence $(F_0,F_1,\dots,F_{t-1},F_0)$ has only one face, a contradiction.

Thus, we may assume $n_1>0$ and $n_2>0$.  In this case, $C$ has at least $n_1\hskip 1pt n_2\ge n-1\ge 4$ crossings, as required for both parts.  
\end{cproof}

%\subsection{Final results}
We have shown in Table~\ref{tbl:count_and_cost} the 
computational cost (accumulated over all cores) for generating
drawings of $K_5$ up to $K_{10}$.   Each drawing $D_{10}$ of $K_{10}$ is used to generate a representative from each equivalence class of drawings in $\widetilde\calD_{11}^{100}$ containing $D_{10}$.   In our
implementation, the representative drawing $D_{11}$ of $K_{11}$ is immediately checked to see if it generates a drawing in
$\calD_{12}^{151}$ having a subdrawing in $\calD_{11}^{104}$.  In particular, the drawing $D_{11}$ is not saved for any later use.

The total time used by Algorithm~\ref{alg:compound} (with in-memory
checking for $K_{12}$) is roughly {20,618
  hours}. With 256 cores available, the average time per core is
less than 4 days.
%<========= ADDED BY SHENGJUN ========

\section{Comments}

The first and third authors are in the process of preparing for publication a computer-free, complete proof that $\crn(K_9)=36$.  (We take the attitude that Guy's proof \cite{guy} is really a computer proof.  For example, he does not demonstrate that there are only three optimal drawings of $K_8$.  Rather, he shows how to do so and did so himself just as a computer would have.)  Some of the considerations required there were useful here; in particular, it was the introduction there of the quantity $\delta$ that led to the surprising insight Lemma \ref{lm:theory}.    

It seems reasonable to assume that Lemma \ref{lm:theory} is just the beginning of progress on $\crn(K_{13})$ and that the work on $K_9$ can be continued to develop more tools to attack both the particular problem of $K_{13}$ as well as Conjecture \ref{co:zaran}.

\section{Acknowledgement}

This work was made possible by the facilities of the Shared Hierarchical Academic Research Computing Network (SHARCNET:www.sharcnet.ca) and Compute/Calcul Canada.   The third author gratefully acknowledges research funding provided by NSERC.

\end{document}